\documentstyle[11pt]{article}

\newtheorem{thm}{Theorem}[section]
\newtheorem{defn}[thm]{Definition}
\newtheorem{prop}[thm]{Proposition}

\newtheorem{cor}[thm]{Corollary}

\def\ba{\begin{array}}
\def\ea{\end{array}}
\def\be{\begin{equation} \label}
\def\ee{\end{equation}}
\def\bea{\begin{eqnarray*}}
\def\eea{\end{eqnarray*}}
\def\beal{\begin{eqnarray} \label}
\def\eeal{\end{eqnarray}}
\def\bit{\begin{itemize}}
\def\eit{\end{itemize}}
\def\rf#1{(\ref{#1})}

\def\a{\alpha}
\def\b{\beta}
\def\g{\gamma}
\def\Ga{\Gamma}
\def\d{\delta}
\def\D{\Delta}

\def\L{\Lambda}
\def\r{\varrho}
\def\s{\sigma}

\def\O{\Omega}

\def\E{{\cal E}}
\def\F{{\cal F}}
\def\G{{\cal G}}
\def\T{{\cal T}}
\def\X{{\cal X}}

\def\Z{{\bf Z}}
\def\Zd{{{\Z}^d}}
\def\R{{\bf R}}
\def\Rd{{{\R}^d}}
\def\GL{{\bf G}}
\def\Lc{{\L^c}}
\def\subsub{\subset\subset}
\def\lr{\leftrightarrow}
\def\llr#1{\stackrel{#1}{\longleftrightarrow}}
\def\p#1{\tilde{#1}}
\def\ph{J}

\def\pp{\mbox{\sf p}}
\def\ss{\mbox{\sf s}}
\def\uu{\mbox{\sf u}}

\begin{document}

\title{{\bf Phase transition and percolation\\ in Gibbsian particle 
models}}
\author{ 
Hans-Otto Georgii\\ {\small\sl Mathematisches Institut der 
Universit\"at M\"unchen}\\ {\small\sl Theresienstr.\ 39, D-80333 
M\"unchen, Germany.} 
}

\date{}
\maketitle

%\begin{quote}
\noindent
We discuss the interrelation between phase transitions 
in interacting lattice or continuum models, and the existence of 
infinite clusters in suitable random-graph models.
In particular, we describe a random-geometric approach to the phase 
transition in the continuum Ising model of two species 
of particles with soft or hard interspecies repulsion. 
We comment also on the related area-interaction process and on 
perfect simulation.
%\end{quote}

\section{Gibbs measures: general principles}

This section contains a brief introduction to the basic physical 
and stochastic ideas leading to the concept of Gibbs measures. The 
principal question is the following:
\begin{quote}
{Which kind of stochastic model is appropriate for the 
description of spatial random phenomena involving a very large number 
of components which are coupled together by an interaction depending 
on their relative position?}
\end{quote}
To find an answer we will start with a spatially discrete 
situation; later we will proceed to the 
continuous case. Consider the phenomenon of ferromagnetism.
\index{ferromagnet}  
A piece of ferromagnetic material like iron or nickel can be imagined 
as 
consisting of many elementary magnets, the so-called spins, which are 
located at the sites of a crystal lattice and have a finite number of 
possible orientations (according to the symmetries of the 
crystal). The essential point is that these spins interact with each 
other in such a way that neighboring spins prefer to be aligned. This 
interaction is responsible for the phenomenon of spontaneous 
magnetization, meaning that at sufficiently low temperatures the 
system can choose between several distinct macrostates in which 
typically all spins have the same orientation. 

How can one find a mathematical model for such a ferromagnet? The 
first fact to observe is that the number of spins is very large.
So, probabilistic experience with the law of large numbers suggests 
to approximate the large finite system by an infinite system in order 
to get clear-cut phenomena. This means that we should assume 
that the underlying crystal lattice is infinite. The simplest case to 
think of is the $d$-dimensional hypercubic lattice $\Zd$. (As the 
case 
$d=1$ is rather trivial, we will always assume that $d\ge2$.)
On the other hand, to keep the model simple it is natural to assume 
that each spin has only finitely many possible orientations. In other 
words, the random spin $\xi_i$ at lattice site $i$  
takes values in a finite state space $S$. The set of all possible 
spin configurations $\xi=(\xi_i)_{i\in\Zd}$ is then the product space 
$\O=S^{\Zd}$. This so-called {\em configuration space} is equipped 
with the Borel $\s$-algebra $\F$ for the natural product topology 
on $\O$. Since the spins are random, we are interested in 
probability measures $P$ on $(\O,\F)$. Such probability spaces are 
known as {\em lattice systems}. For any $\xi\in\O$ and $\L\subset\Zd$ 
we write $\xi_\L=(\xi_i)_{i\in\L}$ for the part of the configuration 
that 
occurs in $\L$. By abuse of notation, we use the same symbol $\xi_\L$ 
for the projection from $\O$ onto $S^\L$.
 
Which kind of probability measure on $(\O,\F)$ can serve as a model 
of 
a ferromagnet? As we have seen above, the essential feature of a 
ferromagnet is the interaction between the spins. We 
are thus interested in probability measures $P$ on $\O$ for which the 
 spin variables $\xi_i$, $i\in\Zd$, are {\em dependent}. 
 A natural way of describing dependencies is to prescribe certain 
 conditional probabilities. This idea, which is familiar from Markov 
 chains, turns out to be suitable also here. Since our parameter set 
$\Zd$ 
 admits no natural linear order, the conditional probabilities can, 
 of course, not lead from a past to a future. Rather we 
{\em prescribe the behavior of a finite set of spins when all other 
spins are fixed}. In other words, we are interested 
 in probability measures $P$ on $(\O,\F)$ having prescribed 
 conditional probabilities 
\be{cond-prob}
G_\L(\xi_\L|\xi_{\L^c})
\ee
for a configuration $\xi_\L\in S^\L$ within a finite set 
$\L\subset\Zd$ given a fixed configuration $\xi_{\L^c}\in S^{\L^c}$ 
off $\L$. In the following we write $\L\subsub\Zd$ when $\L$ is a 
{\em finite\/} subset of $\Zd$. 
The specific form of these conditional distributions does not matter 
at the moment. Two special cases are
\bit
\item {\em the Markovian case: } the conditional distribution 
\rf{cond-prob} only depends on the value of the spins along the {\em 
boundary\/} $\partial\L=\{i\notin\L: |i-j|=1 \mbox{ for some } 
j\in\L\}$ of $\L$, i.e.,
\be{Markov}
G_\L(\xi_\L|\xi_{\L^c})= G_\L(\xi_\L|\xi_{\partial\L})
\ee
(with a slight abuse of notation); $|\cdot|$ stands for the 
Euclidean norm.
\item {\em the Gibbsian case: } the conditional distribution 
\rf{cond-prob} is defined in terms of a {\em Hamilton function\/} 
$H_\L$ by the Boltzmann--Gibbs formula
\be{Gibbs}
G_\L(\xi_\L|\xi_{\L^c})= Z_{\L|\xi_{\L^c}}^{-1} \exp[-H_\L(\xi)] \;,
\ee
where $Z_{\L|\xi_{\L^c}}=\sum_{\xi'\in\O:\xi'\equiv \xi
\mbox{ \scriptsize off }\L}\exp[-H_\L(\xi')]$ is a normalizing 
constant. Physically, $H_\L(\xi)$ describes the 
energy excess of the total configuration $\xi$ over the energy of 
the  outer configuration $\xi_{\L^c}$. (Physicists will miss here 
the factor $\b$, the inverse temperature; we will 
assume that $\b$ is subsumed into $H_\L$ or, equivalently, that 
the units are chosen in such a way that $\b=1$.)
 \eit
 In the following, $G_\L(\,\cdot\,|\xi_{\L^c})$ will be viewed as a 
 probability measure on $\O$ for which the configuration outside $\L$ 
 is almost surely equal to $\xi_\Lc$.
 
 The above idea of prescribing conditional probabilities 
 leads to the following concept  
 introduced in the late 1960's independently by R.L.\ Dobrushin, and 
 O.E.\ Lanford and D.\ Ruelle. 
 \begin{defn}\label{DLR}
 A probability measure $P$ on $(\O,\F)$ is called a Gibbs measure, 
\index{Gibbs measure} or DLR-state, for a 
 family $\GL=(G_\L)_{\L\subsub\Zd}$ of conditional probabilities 
 (satisfying the natural consistency condition) if
 \[
 P(\xi_\L\mbox{ \rm occurs in }\L\,|\, \xi_\Lc \mbox{ \rm occurs off 
}\L)=
 G_\L(\xi_\L\,|\,\xi_{\Lc})  
 \]
 for $P$-almost all $\xi_\Lc$ and all $\L\subsub\Zd$.
 \end{defn}
If $\GL$ is Gibbsian for a Hamiltonian $H$ as in \rf{Gibbs}, each Gibbs 
 measure can be interpreted as an equilibrium state for a physical 
 system with state space $S$ and interaction $H$. This is because 
 the Boltzmann--Gibbs 
 distribution maximizes the entropy when the mean energy is fixed; we 
 will discuss this point in more detail in the continuum setting in
Section \ref{thdyn} below.
 
A general account of the theory of Gibbs measures can be found in the 
monograph Georgii (1988); here we will only present the principal ideas.
In contrast to the situation for Markov chains, Gibbs measures do not 
exist automatically. However, in the present case of a finite state 
space $S$, Gibbs measures do exist whenever $\GL$ is Markovian in the 
sense of \rf{Markov}, or almost Markovian in the sense that the 
conditional probabilities \rf{cond-prob} are continuous functions of 
the outer configuration $\xi_{\L^c}$. In this case one can show that 
any weak 
limit of $G_\L(\,\cdot\,|\xi_\Lc)$ for fixed $\xi\in\O$ as 
$\L\uparrow\Zd$ is 
a Gibbs measure.

The basic observation is that the Gibbs measures for a given 
consistent family $\GL$ of conditional probabilities form a convex 
set $\G$. Therefore one is interested in its extremal points. These 
can be characterized as follows.
\begin{thm} Let $\T=\bigcap\s(\xi_\Lc:\L\subsub\Zd)$ the tail 
$\s$-algebra, i.e., the $\s$-algebra of all macroscopic events not 
depending on the values of any finite set of spins. Then the 
following statements hold:

{\rm (a)} A Gibbs measure $P\in\G$ is extremal in $\G$ if and only if $P$ 
is trivial on $\T$, i.e., if and only if any tail measurable real 
function is $P$-almost surely constant.

{\rm (b)} Any two distinct extremal Gibbs measure are mutually singular on 
$\T$.

{\rm (c)} Any non-extremal Gibbs measure is the barycenter of a unique 
probability weight on the set of extremal Gibbs measures.
\end{thm}
A proof can be found in Georgii (1988), Theorems (7.7) and (7.26).
Statement (a) means that the extremal Gibbs measures are {\em 
macroscopically deterministic\/}: on the macroscopic level all 
randomness disappears, and an experimenter will get non-fluctuating 
measurements of macroscopic quantities like magnetization 
or energy per lattice site. Statement (b) asserts that {\em distinct 
extremal Gibbs measures show different macroscopic behavior}. So, 
they can be distinguished by looking at typical realizations of the 
spin configuration through macroscopic glasses. Finally, statement 
(c) implies that any realization which is typical for a non-extremal 
Gibbs measure is in fact typical for a suitable extremal Gibbs measure. 
In physical terms: {\em any configuration which can be seen in nature 
is governed 
by an extremal Gibbs measure}, and the non-extremal Gibbs measures 
can only be interpreted in a Bayesian way as measures describing the 
uncertainty of the experimenter. These observations lead us to the 
following definition.
\begin{defn}
Any extremal Gibbs measure is called a phase of the corresponding 
physical system. If distinct phases exist, one says that a phase 
transition occurs. \index{phase transition}
\end{defn} 
So, in terms of this definition the existence of phase transition is 
equivalent to the non-uniqueness of the Gibbs measure. In the light 
of 
the preceding theorem, this corresponds to a ``macroscopic 
ambivalence'' 
of the system's behavior. We should add that not all critical 
phenomena in nature can be described in this way: even when the Gibbs 
measure is unique it may occur that it changes its qualitative 
behavior when some parameters of the interaction are changed. 
However, we will not discuss these possibilities here and stick to 
the definition above for definiteness.
In this contribution we will ask:
\begin{quote}
What are the driving forces giving rise to a phase transition?
Is there any mechanism relating microscopic and macroscopic behavior 
of spins?
\end{quote}
As we will see, in a number of cases one can give the following geometric 
answer: 
\begin{quote}
One such mechanism is the formation of infinite clusters 
in suitably defined random graphs. Such infinite clusters serve as a link 
between the local and global behavior of spins, and make visible 
how the individual spins unite to form a 
specific collective behavior.
\end{quote}
In the next section we will discuss two lattice models for which 
this answer 
is correct. In Section 4 we will show that a similar answer 
can also be given in a continuum set-up. A useful 
technical tool is the stochastic comparison of probability measures.
 
Suppose the state space $S$ is a subset of $\R$ and thus linearly 
ordered. Then the configuration space $\O$ has a natural partial 
order, 
and we can speak of increasing real functions. Let $P,P'$ be two 
probability measures on $\O$. We say that $P$ is {\em stochastically 
smaller} than $P'$, and write $P\preceq P'$, if $\int f\,dP\le \int 
f\, dP'$ \index{stochastic order}
for all local increasing functions (or, equivalently, for all 
measurable bounded increasing functions) $f$ on $\O$. A sufficient 
condition for stochastic monotonicity is given in the proposition 
below. Although this 
condition refers to the case of finite products (for which 
stochastic monotonicity is similarly defined), it is also useful in 
the case of infinite product spaces. This is because (by the very 
definition) the relation $\preceq$ is preserved under weak limits.
\begin{prop}\label{Holl}
{\bf(Holley's inequality) }Let $S$ be a finite subset of $\R$, $\L$ 
a finite index set, and $P,P'$ two probability measures on the finite 
product space $S^\L$ giving positive weight to each element of $S^\L$.
Suppose the single-site conditional probabilities at any $i\in\L$ 
satisfy
\[
P(\,\cdot\,|\xi_{\L\setminus\{i\}}\mbox{ \rm occurs off $i$} ) 
\preceq 
P'(\,\cdot\,|\xi'_{\L\setminus\{i\}}\mbox{ \rm occurs off $i$})
\quad\mbox{ whenever $\xi\le \xi'$.}
\]
Then $P\preceq P'$. If this condition holds with $P'=P$ then $P$ has 
positive correlations in the sense that any two bounded increasing 
functions are positively correlated.
\end{prop} 
For a proof (and a slight extension) we refer to Theorems 4.8 and 
4.11 of Georgii, H\"aggstr\"om and Maes (1999).

\section{Phase transition and percolation: two lattice models}

To provide the necessary background for our results on continuum 
particle systems let us still stick to the lattice case. We will 
consider two classical models which allow an understanding of phase 
transition in random-geometric terms. Many further examples for the 
relation between random geometry and 
phase transition can be found in Georgii, H\"aggstr\"om and Maes (1999).

Let us start recalling some basic facts
on {\em Bernoulli percolation} on $\Zd$ for $d\ge2$.
\index{percolation!Bernoulli site bond|(}
Consider $\Zd$ as a graph with vertex 
set $\Zd$ and edge set $E(\Zd)=\{e=\{i,j\}\subset\Zd:|i-j|=1\}$.
We fix two parameters $0\le p_s,p_b\le1$, the site and bond 
probabilities, and construct a random subgraph $\Ga=(X,E)$ of 
$(\Zd,E(\Zd))$ by setting
\[
X=\{i\in\Zd: \xi_i=1\}\,, \quad E=\{e\in E(X): \eta_e=1\}\,,
\]
where $E(X)=\{e\in E(\Zd): e\subset X\}$ is the set of all edges 
between the sites of $X$, and $\xi_i$, $i\in\Zd$, and $\eta_e$, $e\in 
E(\Zd)$, are independent Bernoulli variables satisfying 
$P(\xi_i=1)=p_s$, $P(\eta_e=1)=p_b$.
This construction is called the {\em Bernoulli mixed site-bond 
percolation model\/}; setting $p_b=1$ we obtain pure site percolation, and 
the case $p_s=1$ corresponds to pure bond percolation.

Let $\{0\lr\infty\}$ denote the event that $\Ga$ contains an infinite 
path starting from $0$, and 
$$\theta(p_s,p_b;\Zd)=\mbox{Prob}(0\lr\infty)$$
be its probability. By Kolmogorov's zero-one law, we have 
$\theta(p_s,p_b;\Zd)>0$ if and only if $\Ga$ contains an infinite cluster 
with probability $1$. 
In this case one says that percolation occurs.
The following proposition asserts that this happens in a non-trivial 
region of the parameter square, which is separated by the so-called 
critical line from the region where all clusters of $\Ga$ are almost 
surely finite.
The change of behavior at the critical line is the simplest 
example of a critical phenomenon. 
\begin{prop}\label{perc} 
The function $\theta(p_s,p_b;\Zd)$ is increasing in $p_s$, $p_b$ and 
$d$. Moreover,
$\theta(p_s,p_b;\Zd)=0$ when $p_sp_b$ is small enough, while 
$\theta(p_s,p_b;\Zd)>0$ when $d\ge2$ and $p_sp_b$ is sufficiently close to $1$.
\end{prop}
{\em Sketch proof: }The monotonicity in $p_s$ and $p_b$ follows 
from Proposition \ref{Holl}, and the one in $d$ from an obvious 
embedding argument. To show that $\theta=0$ when $p_sp_b$ is 
small, we note that the expected number of neighbors in $\Ga$ of a 
given lattice site is $2dp_sp_b$. Comparison with a branching process 
thus shows that $\theta=0$ when $2dp_sp_b<1$.

Next, let $d=2$ and suppose $0\in X$ but the cluster $C_0$ of $\Ga$ 
containing $0$ is finite. Consider $\partial_{ext}C_0$, the part of 
$\partial C_0$ belonging to the infinite component of $C_0^c$. For 
each site $i\in \partial_{ext}C_0$, either this site or all bonds 
leading from $i$ to $C_0$ do not belong to $\Ga$. This occurs with 
probability at most $1-p_sp_b$. So, the probability that 
$\partial_{ext}C_0$ has a fixed location is at most 
$(1-p_sp_b)^\ell$ with $\ell=\#\partial_{ext}C_0$.
Counting all possibilities for this location one finds that 
$1-\theta<1$ when $1-p_sp_b$ is small enough. 
By the monotonicity in 
$d$, the same holds a fortiori in higher dimensions. $\Box$
\index{percolation!Bernoulli site bond|)}
 
\medskip\noindent
The above proposition is all what we need here on Bernoulli 
percolation; an excellent source for a wealth of further rigorous 
results is the book of Grimmett (1999).

We now ask for the role of percolation for Gibbs measures, and in 
particular for the existence of phase transitions. Of course, 
in contrast to the above this will involve {\em dependent\/}, i.e., 
non-Bernoulli percolation. We consider here two specific examples. 
 In  
these examples, the family $\GL$ of conditional probabilities 
will be Gibbsian for a nearest-neighbor 
interaction; this means that both \rf{Markov} and \rf{Gibbs} are 
valid.

\subsection{The Ising model}
\index{Ising model!lattice|(}

This is by far the most famous model of Statistical Mechanics, named 
after E.\ Ising who studied this model in the early 1920s in his 
thesis suggested by W.\ Lenz. It is the simplest model of a 
ferromagnet in equilibrium. One assumes that the spins have only two 
possible orientations, and therefore defines $S=\{-1,1\}$. The 
family $\GL$ is defined by \rf{Gibbs} with
\be{Ising}
H_\L(\xi)= J \sum_{\{i,j\}\cap\L\ne\emptyset:|i-j|=1} 1_{\{\xi_i\ne 
\xi_j\}}\;,
\ee
where $J>0$ is a coupling constant which is inversely proportional to 
the absolute temperature. This means that neighboring spins of 
different sign have to pay an energy cost $J$. There exist two 
configurations of minimal energy, so-called ground 
states, namely the configuration `{$+$}' which is 
identically equal 
to $+1$, and the configuration `{$-$}' identically equal 
to $-1$. The behavior of the model is governed by these two ground 
states. To see this we begin with some useful consequences of the 
ferromagnetic character of the interaction. First, it is intuitively 
obvious that the measures $G_\L^+=G_\L(\,\cdot\,|+)$ 
decrease stochastically when $\L$ increases (since then the 
effect of the {$+$} boundary decreases). This follows easily 
from Holley's inequality, Proposition \ref{Holl}. Since the local 
increasing functions are a 
convergence determining class, it follows that the 
weak infinite-volume limit 
$P^+=\lim_{\L\uparrow\Zd}G_\L^+$ exists. 
Likewise, the weak limit 
$P^-=\lim_{\L\uparrow\Zd}G_\L(\,\cdot\,|-)$ exists (and 
by symmetry is the image of $P^+$ under simultaneous spin flip). 
These limits are Gibbs measures and invariant under translations. 
Holley's inequality also implies that they are 
stochastically maximal resp.\ minimal in $\G$, and 
in particular extremal. This gives us the following criterion for 
phase transition in the Ising model.
\begin{prop}\label{phtr-criterium}
For the Ising model on $\Zd$ with any coupling constant $J>0$ we have
$\#\G>1$ if and only if $P^-\ne P^+$ if and only if 
$\int \xi_0\,dP^+>0$.
\end{prop}
The last equivalence follows from the relation $P^-\preceq P^+$, 
the translation invariance of these Gibbs measures, and the spin-flip 
symmetry. A detailed proof of the proposition and the previous 
statements can be found in Section 4.3 of Georgii, H\"aggstr\"om and Maes 
(1999).

How can we use this criterion? This is where random 
geometry enters the scenery. The key is the following geometric 
construction tracing back to Fortuin and Kasteleyn (1972) 
and in this form to Edwards and Sokal (1988). It is called 
the {\em random-cluster representation of the Ising model}.
\index{random-cluster representation!lattice Ising model|(}

Let $\E_\L^+=\{E\subset E(\Zd):E\supset E(\Lc)\}$ be the set of all edge 
configurations in $\Zd$ which include all edges outside $\L$, 
and define a probability measure $\phi_\L$ on $\E_\L^+$ by setting
\be{RC}
\phi_\L(E) = Z_\L^{-1}\; 2^{k(E)}\, 
p^{\#E \setminus E(\Lc)}\,(1-p)^{\# E(\Zd)\setminus E} \quad\mbox{ 
when $E\supset E(\Lc)$,}
\ee
where $p=1-e^{-J}$, $k(E)$ is the number of clusters of the graph
$(\Zd,E)$, and $Z_\L$ is a normalizing constant. 
$\phi_\L$ is called the {\em random-cluster distribution in $\L$ 
with wired boundary condition}. This measure turns out to be
related to $G_\L^+=G_\L(\,\cdot\,|+)$. It will be 
convenient to identify a configuration $\xi\in\O$ with the pair 
$(X^+,X^-)$, where $X^+$ and $X^-$ are the sets of all lattice sites 
$i$ for which $\xi_i=+1$ resp.\ $-1$. 
\begin{prop}\label{prop:RC} 
For any hypercube $\L$ in $\Zd$ there exists the following 
correspondence 
between the the Gibbs distribution $G_\L^+$ for the Ising model
and the random-cluster distribution $\phi_\L$ in \rf{RC}.

$(\,G_\L^+\leadsto\phi_\L\,) $ Take a spin configuration 
$\xi=(X^+,X^-)\in\O$ with distribution $G_\L^+$, and define an 
edge configuration $E\in\E_\L^+$ as follows:
Independently for all $e\in E(\Zd)$  
let $e\in E$ with probability 
\[
p_\L(e)=\left\{\ba{cl}
1-e^{-J}&\mbox{\rm if $e\subset X^+$ or $e\subset X^-$, and 
$e\cap\L\ne\emptyset$,}\\ 
1&\mbox{\rm if $e\subset\Lc$,}\\
0&\mbox{\rm otherwise,} \ea\right.
\]
and $e\notin E$ otherwise. Then $E$ has distribution $\phi_\L$.

$(\,\phi_\L\leadsto G_\L^+\,) $  Pick an edge configuration 
$E\in\E_\L^+$ according to $\phi_\L$, and define a spin configuration 
$\xi=(X^+,X^-)\in\O$ as follows:
For each finite cluster $C$ of $(\Zd,E)$ let $C\subset X^+$
or $C\subset X^-$ according to independent flips of a fair coin; 
the unique infinite cluster of $(\Zd,E)$ containing $\Lc$ is included 
into $X^+$. Then $\xi$ has distribution $G_\L^+$.
\end{prop}
{\em Proof. }A joint description of a spin configuration $\xi\in\O$ 
with 
distribution $G_\L^+$ and an edge configuration $E\in\E_\L^+$ with 
distribution $\phi_\L$ can be obtained as follows.
Any edge $e\in E(\Zd)$ is independently included into $E$ 
with probability $p=1-e^{-J}$ resp.\  $1$ according to whether 
$e\cap\L\ne\emptyset$ or not; each spin in $\L$ is 
equal to $+1$ or $-1$ according to independent flips of a fair coin;
the spins off $\L$ are set equal to $+1$. The measure $P$ thus 
described is then conditioned on the event $A$ that no spins of 
different sign are connected by an edge. 
Relative to $P(\,\cdot\,|A)$, $\xi$ has distribution 
$G_\L^+$ and $E$ has distribution $\phi_\L$.
This is because $\exp[-H_\L(\xi)]$ is equal to the conditional 
$P$-probability of $A$ given $\xi$, and
$2^{k(E)-1}$ is proportional to the conditional $P$-probability of $A$ 
given $E$. Now, it is easy to see that the two constructions in the 
proposition simply correspond to the conditional distributions of $E$ 
given $\xi$ resp.\ of $\xi$ given $E$ relative to $P(\,\cdot\,|A)$.
$\Box$

\medskip\noindent
Intuitively, the edges in the random-cluster representation indicate 
which pairs of spins ``realize'' their interaction, in that they 
decide to take the same orientation to avoid the dealignment costs.
On the one hand, this representation is the basis of an efficient 
simulation procedure, the algorithm of Swendsen--Wang (1987),
\index{Swendsen--Wang algorithm}
which together with its continuous counterpart will be discussed at the 
end of Section \ref{simul}. On the other hand, it is the key for a
geometric approach to the phenomenon of phase transition, as we will
now show.

The construction $\phi_\L\leadsto G_\L^+$ implies that,
for $0\in\L$, the conditional expectation 
of $\xi_0$ given $E$ is $1$ when $0$ is connected to $\partial\L$ by 
edges in $E$, and $0$ otherwise. Hence
\be{magn-percol1}
\int \xi_0\,dP_\L^+=\phi_\L(0\lr\partial\L)\;.
\ee
By Proposition \ref{Holl}, the measures $\phi_\L$ decrease 
stochastically 
when $\L$ increases, so that the infinite-volume 
random-cluster distribution $\phi=\lim_{\L\uparrow\Zd}\phi_\L$ 
exists. 
Letting $\L\uparrow\Zd$ in \rf{magn-percol1} we thus find that
$
\int \xi_0\,dP^+=\phi(0\lr\infty)\,.
$
\index{random-cluster representation!lattice Ising model|)}
Combining this with Proposition \ref{phtr-criterium} we obtain the
first statement of the following theorem, the 
equivalence of percolation and phase transition. This gives us 
detailed 
information on the existence of phase transition. 
\begin{thm}\label{Ising-phtr}
Consider the Ising model on $\Zd$ with 
Hamiltonian \rf{Ising} for any coupling 
constant $J>0$. Then $\#\G>1$ if and only if $\phi(0\lr\infty)>0$.
Consequently, there exists  a coupling threshold $0<J_c<\infty$ 
(corresponding to a critical inverse temperature) such that $\#\G=1$ 
when $J<J_c$ and $\#\G>1$ when $J>J_c$.
\end{thm} \index{phase transition!lattice Ising model}
{\em Sketch proof. }It only remains to show the second statement. 
This follows from Holley's inequality, Proposition \ref{Holl}. 
First, this inequality implies 
that $\phi_\L$ is stochastically increasing in the parameter 
$p=1-e^{-J}$. Hence $\phi(0\lr\infty)$ is an increasing function of 
$p$. Moreover, one finds that $\phi$ is stochastically dominated by 
the Bernoulli bond percolation measure, whence $\phi(0\lr\infty)=0$ 
when $p$ is so small that $\theta(1,p;\Zd)=0$. 
Finally, $\phi$ is stochastically larger than the Bernoulli 
bond percolation measure with parameter $\tilde p=p/(p+2(1-p))$. 
Hence $\phi(0\lr\infty)>0$ when $p$ is so large that 
$\theta(1,\tilde p;\Zd)>0$. Details of this computation 
can be found in Section 6 of Georgii, H\"aggstr\"om and Maes (1999), 
which deals in fact with the 
extension of these results to the Potts model in which each spin has 
$q$ different values. $\Box$

\medskip\noindent
One may ask whether the connection between percolation and phase 
transition can be seen more directly in the behavior of spins. The 
following corollary gives an answer to this question.
Let $\{0\llr{+}\infty\}$ denote the event 
that $0$ belongs to an infinite cluster of the graph $(X^+,E(X^+))$ 
induced by the set of plus-spins.
\begin{cor}\label{Ising-plus-perc} 
For the Ising model on $\Zd$ with arbitrary coupling 
constant $J>0$ we have $P^+(0\llr{+}\infty)>0$ whenever
$\#\G>1$. The converse holds only when $d=2$.
\end{cor} \index{percolation!lattice Ising model}
The first part follows readily from the construction 
in Proposition \ref{prop:RC} which shows
that $P^+(0\llr{+}\infty)\ge \phi(0\lr\infty)$. For its second part 
see Georgii, H\"aggstr\"om and Maes (1999).
Pursuing the idea of plus-percolation further one can obtain the
following result independently obtained in the late 
1970s by Aizenman and Higuchi on the basis of previous work of L.\ 
Russo; a simpler proof has recently been given by Georgii and Higuchi 
(1999).
\begin{thm}\label{2d} 
For the Ising model on $\Z^2$ with $J>J_c$, there exist no other phases 
than $P^+$ and $P^-$.
\end{thm}
A celebrated result of Dobrushin asserts that in three or more 
dimensions 
there exist non-translation invariant Gibbs measures which look like 
$P^+$ in one half-space and like $P^-$ in the other half-space.
\index{Ising model!lattice|)}

\subsection{The Widom--Rowlinson lattice gas}
\index{Widom--Rowlinson!lattice model|(}

The Widom--Rowlinson lattice gas is a discrete analog of a continuum 
model to be considered in Section 4. It describes the random 
configurations of particles of two different types, plus or minus, which 
can only sit at the sites of the lattice $\Zd$. Multiple occupations 
are excluded. So, at each site $i$ of the 
lattice there are three possibilities: either $i$  
is occupied by a plus-particle, or by a minus-particle, or $i$ is 
empty. 
The configuration space is thus $\O=S^\Zd$ with $S=\{-1,0,1\}$.
The basic assumption is that there is a hard-core repulsion between 
plus- and minus-particles, which means that particles of 
distinct type are not allowed to sit next to each other. 
In addition, 
there exists a chemical ``activity'' $z>0$ which governs the 
overall-density of particles. The Hamiltonian thus takes the form
\be{WR}
H_\L(\xi)= \sum_{\{i,j\}\cap\L\ne\emptyset:|i-j|=1} U(\xi_i,\xi_j)
-\log z\;\sum_{i\in\L} |\xi_i|\;,
\ee
where $U(\xi_i,\xi_j)=\infty$ if $\xi_i \xi_j=-1$, and 
$U(\xi_i,\xi_j)=0$ otherwise. The associated family $\GL$ of 
conditional probabilities is again given by \rf{Gibbs}. Just as in the 
Ising model, for $z>1$ there exist two distinguished configurations 
of minimal energy, namely the constant configurations 
`{$+$}' and `{$-$}' for which all sites are occupied 
by particles of the same type. Moreover, one can again apply
Holley's inequality to show that the Gibbs distributions
$G_\L^+=G_\L(\,\cdot\,|+)$ and
$G_\L^-=G_\L(\,\cdot\,|-)$ converge to translation 
invariant limits $P^+,P^-\in\G$ which are 
stochastically maximal resp.\ minimal in $\G$, and 
therefore extremal. This implies that
Proposition \ref{phtr-criterium} holds verbatim also in the present 
case.

Is there also a geometric representation of the model, just as for 
the Ising model? The answer is yes, with interesting analogies and 
differences. There exists again a
\index{random-cluster representation!Widom--Rowlinson lattice model|(}
random-cluster distribution with an appearance very similar to 
\rf{RC}, but this involves site percolation rather than bond 
percolation. Namely, consider the probability measure 
$\psi_\L$
on the set $\X_\L^+=\{Y\subset\Zd:Y\supset \Lc)\}$  which is given by
\be{WR-RC}
\psi_\L(Y) = Z_\L^{-1}\; 2^{k(Y)}\, 
p^{\#Y \setminus \Lc}\,(1-p)^{\# \Zd\setminus Y} \quad\mbox{ 
for $Y\supset \Lc$\,;}
\ee
here $p=\frac{z}{1+z}$, $k(Y)$ is the number of clusters of the graph
$(Y,E(Y))$, and $Z_\L$ is a normalizing constant. 
$\psi_\L$ is called the 
{\em site random-cluster distribution in $\L$ with parameter $p$ and 
wired boundary condition}. We will again identify a configuration 
$\xi\in\O$ with a pair $(X^+,X^-)$, 
where $X^+$ and $X^-$ are the sets of all lattice points 
$i$ such that $\xi_i=+1$ resp.\ $-1$; thus $X^+\cup X^-$ is the set of 
occupied sites, and $\xi\equiv 0$ on its complement.
Here is the {\em random-cluster representation of the 
Widom--Rowlinson 
lattice gas} which is analogous to Proposition \ref{prop:RC}.
\begin{prop}\label{prop:WR-RC} 
For any hypercube $\L$ in $\Zd$ there exists the following 
correspondence between the the Gibbs distribution $G_\L^+$ 
for the Widom--Rowlinson model 
and the site random-cluster distribution 
$\psi_\L$ in \rf{WR-RC}.

$(\,G_\L^+\leadsto\psi_\L\,) $ For a random spin configuration 
$\xi=(X^+,X^-)\in\O$ with distribution $G_\L^+$, the random set
$Y=X^+\cup X^-$ has distribution $\psi_\L$.

$(\,\psi_\L\leadsto G_\L^+\,) $ Pick a random set 
$Y\in\X_\L^+$ according to $\psi_\L$, and define a spin configuration 
$\xi=(X^+,X^-)\in\O$  with $X^+\cup X^-=Y$ as follows:
For each finite cluster $C$ of $(Y,E(Y))$ let $C\subset X^+$
or $C\subset X^-$ according to independent flips of a fair coin; 
the unique infinite cluster of $(Y,E(Y))$ containing $\Lc$ is 
included into $X^+$. Then $\xi$ has distribution $G_\L^+$.
\end{prop}
The random-cluster representation of the Widom-Rowlinson model is 
simpler than that of the Ising model because the randomness 
involves only the sites of the lattice. (This is a consequence of the 
hard-core interaction; in the case of a soft repulsion the situation 
would be different, as we will see in the continuum setting in 
Section 4.) On the other hand, there is a serious drawback of
the site random-cluster distribution $\psi_\L$:
it does not satisfy the conditions of Proposition \ref{Holl} for 
positive correlations. This is because the conditional probabilities 
in \rf{psi-cond-prob} below are not increasing in $Y$. So, we 
still have a counterpart to \rf{magn-percol1}, viz.\
\be{WR-percol}
\int \xi_0\,dG_\L^+=\psi_\L(0\lr\partial\L)\;,
\ee
but we do not know if the measures $\psi_\L$ are 
stochastically increasing in $z$ and stochastically decreasing in $\L$.
\index{random-cluster representation!Widom--Rowlinson lattice model|)} 
So we obtain a somewhat weaker theorem. 
\begin{thm}\label{WR-phtr}
Consider the Widom--Rowlinson model on $\Zd$, $d\ge2$, 
defined by \rf{WR} with activity 
$z>0$. Then $\#\G>1$ if and only if 
$\lim_{\L\uparrow\Zd}\psi_\L(0\lr\partial\L)>0$. In particular, we 
have $\#\G=1$ when $z$ is sufficiently small, and $\#\G>1$ when 
$z$ is large enough.
\end{thm} \index{phase transition!Widom--Rowlinson lattice model}
{\em Sketch proof. }The first statement follows immediately from
\rf{WR-percol} and the analog of Proposition \ref{phtr-criterium}.
To prove the second statement we note that $\psi_\L$ has single-site 
conditional probabilities of the form
\be{psi-cond-prob}
\psi_\L(Y\ni i\,|\,Y\setminus\{i\})=
{p}\Big/[{p+(1-p)\,2^{\kappa(i,Y)-1}}]\;,
\ee 
where $\kappa(i,Y)$ is the number of clusters of 
$Y\setminus\{i\}$ that intersect a neighbor of $i$. Since
$0\le \kappa(i,Y)\le 2d$, it follows from Proposition \ref{Holl} that 
$\psi_\L$ is stochastically dominated by the site-Bernoulli measure 
with 
parameter $p^*=p/(p+(1-p)2^{-1})$, and dominates the site-Bernoulli 
measure with parameter $p_*=p/(p+(1-p)2^{2d-1})$. Combining 
this with Proposition \ref{perc} we thus find that 
$\#\G=1$ when $z$ is so small that $\theta(p^*,1;\Zd)=0$, and $\#\G>1$ 
when $z$ is so large that $\theta(p_*,1;\Zd)>0$. $\Box$

\medskip\noindent
Let us note that, in contrast to Theorem \ref{Ising-phtr}, the 
preceding result does not extend to the case when there are more than 
two different types of particles; this is related to the lack of 
stochastic monotonicity in this model. However, due to its 
simpler random-cluster representation the 
Widom--Rowlinson model has one advantage over the Ising model, in 
that it satisfies a much stronger counterpart to 
Corollary \ref{Ising-plus-perc}. In analogy to the notation there we 
write $\{0\llr{+}\infty\}$ for the event 
that the origin belongs to an infinite cluster of plus-particles.
\begin{cor}\label{WR-plus-perc} 
For the Widom--Rowlinson lattice gas on $\Zd$ for 
arbitrary dimension $d\ge2$ and with any 
activity $z>0$, we have $\#\G>1$ if and only if 
$P^+(0\llr{+}\infty)>0$.
\end{cor}
\index{percolation!Widom--Rowlinson lattice model}
{\em Sketch proof. }The construction in Proposition \ref{prop:WR-RC} 
readily implies that 
\[
\psi_\L(0\lr\partial\L)=G_\L^+(0\llr{+}\partial\L)\;.
\]
Combining this with \rf{WR-percol} and letting $\L\uparrow\Zd$ one 
obtains the result. $\Box$

\medskip\noindent The above equivalence of phase transition and 
percolation even holds when $\Zd$ is replaced by an arbitrary graph. 

As noticed before Theorem \ref{WR-phtr}, we have no stochastic 
monotonicity in the activity $z$, and therefore no activity threshold 
for the existence of a phase transition. We are thus led to ask if, 
at least, the particle density is an increasing function of $z$.
\index{particle density!monotonicity}
This can be deduced from general thermodynamic principles relying 
on convexity of thermodynamic functions rather than stochastic 
monotonicity. This will be described in Section \ref{thdyn} 
in the continuum set-up.
\index{Widom--Rowlinson!lattice model|)}

\section{Continuum percolation}
\index{percolation!Poisson random-edge|(}

In the rest of this contribution we will show that quite a lot of
the preceding results and techniques carry over to models of point 
particles in Euclidean space. In this section we deal with a simple 
model of continuum percolation. Roughly speaking, this model 
consists of Poisson points which are connected by Bernoulli edges.
To be precise, let $\X$ denote the set of all locally finite subsets 
$X$ of $\R^d$. $\X$ is the set of all point configurations in $\R^d$,
and is equipped with the usual $\s$-algebra generated by the 
counting variables $X\to\#X_\L$ for bounded Borel sets 
$\L\subset\R^d$;
here we use the abbreviation $X_\L=X\cap\L$. 
Next, let $\E$ be the set of all locally finite 
subsets of $E(\R^d)=\{\{x,y\}\subset\R^d:x\ne y\}$. $\E$ is the set 
of all possible edge configurations and is equipped 
with an analogous $\s$-algebra. For $X\in\X$ let $E(X)=\{e\in 
E(\R^d):e\subset X\}$ the set of all possible edges between the 
points 
of $X$, and $\E_X=\{E\in\E:E\subset E(X)\}$ the set of edge 
configurations between the points of $X$.
We construct a random graph $\Ga=(X,E)$ in $\R^d$ as follows. 
\bit
\item Pick a random point configuration $X\in\X$ according to the 
Poisson point process $\pi^z$ on $\R^d$ with intensity $z>0$.
\item For given $X\in\X$, pick a random edge configuration $E\in\E_X$
according to the Bernoulli measure $\mu_X^p$ on $\E_X$ for which the 
events $\{E\ni e\}$, $e=\{x,y\}\in E(X)$, are independent with 
probability $\mu_X^p(E\ni e)=p(x-y)$; here $p:\R^d\to[0,1]$ 
is a given even measurable function.
\eit
The distribution of our random graph $\Ga$ is thus determined by the 
probability measure
\be{rand-edge}
P^{z,p}(dX,dE)=\pi^z(dX)\,\mu_X^p(dE)
\ee
on $\X\times\E$. It is called the {\em Poisson 
random-edge model}, or {\em Poisson random-connection model}, and has 
been introduced and studied first by M.\ Penrose (1991); a 
detailed account of its properties is given in Meester and Roy (1996).

A special case of particular interest is when $p(x-y)$ is equal to 
$1$ when $|x-y|\le 2r$ for some $r>0$, and $0$ 
otherwise. This means that any two points $x$ and $y$ are connected 
by an edge if and only if the balls $B_r(x)$ and $B_r(x)$ 
with radius $r$ and
center $x$ resp.\ $y$ overlap. The connectivity properties of the 
corresponding Poisson random-edge model are thus the same as those of 
the random set $\Xi=\bigcup_{x\in X}B_r(x)$, for random $X$ with 
distribution $\pi^z$. This special case is therefore called the 
{\em Boolean model}, or the {\em Poisson blob model}.
\index{percolation!Boolean}

Returning to the general case, we consider the {\em percolation 
probability\/} of a typical point. Writing $x\lr\infty$ when 
$x$ belongs to an infinite cluster of $\Ga=(X,E)$, this is given by the 
expression 
\be{cont-theta}
\theta(z,p;\Rd)=\int \frac{\#\{x\in X_\L: 
x\lr\infty\}}{|\L|}\;
P^{z,p}(dX,dE)
\ee
for an arbitrary bounded box $\L$ with volume $|\L|$.
By translation invariance, $\theta(z,p;\Rd)$ does not depend on $\L$. 
In fact, in 
terms of the Palm measure $\hat\pi^z$ of $\pi^z$ and the associated 
measure $\hat P^{z,p}(dX,dE)=\hat\pi^z(dX)\,\mu_X(dE)$ we can write
\[
\theta(z,p;\Rd)=\hat P^{z,p}(0\lr\infty)\;.
\]
The following result of M.\ Penrose (1991) is the continuum 
analog of Proposition \ref{perc}.
\begin{thm}\label{cont-perc} 
$\theta(z,p;\Rd)$ is an increasing function of the intensity
$z$ and the edge probability function $p(\cdot)$. Moreover, 
$\theta(z,p;\Rd)=0$
when $\,z\int p(x)\,dx\,$ is sufficiently small, while
$\theta(z,p;\Rd)>0$ when $\,z\int p(x)\,dx\,$ is large enough.
\end{thm}
{\em Sketch proof: }The monotonicity follows from an obvious stochastic 
comparison argument. Since $z\int p(x)\,dx$ is the 
expected number of edges emanating from a given point, a 
branching argument shows that $\theta(z,p;\Rd)=0$ when $z\,\int 
p(x)\,dx<1$. It remains to show that $\theta(z,p;\Rd)>0$ when 
$z\,\int p(x)\,dx$ is 
large enough. By scaling we can assume that $\int 
p(x)\,dx=1$. For simplicity we will in fact suppose that $p$ is bounded 
away from $0$ in a neighborhood of the origin, i.e.,
$p(x-y)\ge\d>0$ whenever $|x-y|\le 2r$. (The following is a special 
case of an argument of Georgii and H\"aggstr\"om (1996).)

We divide the space $\R^d$ into cubic cells $\D(i)$, $i\in\Zd$, 
with diameter at most $r$. We also pick a sufficiently large number 
$n$ and introduce the following two concepts.
\bit
\item Call a cell $\D(i)$ {\em good\/} if it contains at least $n$ 
points which form a connected set relative to the edges of $\Ga$ 
in between them. 
This event does not depend on the configurations in all other cells 
and has probability at least
\[
\pi^z(N_i\ge n)\;[1-(n-1)(1-\d^2)^{n-2}]\equiv p_s\;;
\]
here, $N_i$ is the random number of points in cell $\D(i)$, and the 
second term in the square bracket is an estimate for the probability 
that one 
of the $n$ points is not connected to the first point by a sequence 
of two edges. The essential fact is that $p_s$ is arbitrarily close 
to $1$ when $n$ and $z$ are large enough.
\item Call two adjacent cells $\D(i),\D(j)$ {\em linked\/} if there 
exists an edge from some point in $\D(i)$ to some point in $\D(j)$. 
Conditionally on the event that $\D(i)$ and $\D(j)$ are good, 
this has probability at least $1-(1-\d)^{n^2}=p_b$, 
which is also close to $1$ when $n$ is large enough.
\eit
Now the point is the following: whenever there exists an infinite 
cluster of linked good cells (i.e., an infinite cluster in the 
countable graph with vertices at the good cells and with edges 
between pairs of linked cells) then there exists an infinite cluster 
in the original Poisson random-edge model. Hence $\theta(z,p;\Rd)\ge 
\frac nv\,\theta(p_s,p_b;\Zd)$, where $v$ is the cell volume and 
$\theta(p_s,p_b;\Zd)$ is as in Proposition 
\ref{perc}. Hence $\theta(z,p;\Rd)>0$ when $z$ is large enough. $\Box$
\index{percolation!Poisson random-edge|)}

\medskip\noindent
How can one extend a percolation result as above from the Poisson 
case to point processes with spatial dependencies? Just as in the 
lattice gas, one can take advantage of stochastic comparison 
techniques. 
To this end we need a continuum analog of Holley's theorem.

A simple point process $P$ on a bounded Borel subset $\L$ of $\Rd$ 
(i.e., a probability measure on $\X_\L=\{X\in\X:X\subset\L\}$) is 
said to have Papangelou (conditional) intensity $\g:\L\times\X_\L\to 
[0,\infty[$ if $P$ satisfies the identity
\be{Papangelou}
\int P(dX) \sum_{x\in X} f(x,X\setminus\{x\})=
\int dx\int P(dX)\, \g(x|X)\, f(x,X)
\ee
for any measurable function $f:\L\times\X_\L\to [0,\infty[$. (This is 
a non-stationary analog of the Georgii--Nguyen--Zessin equality 
discussed in the contribution of D.\ Stoyan to this volume.) 
This equation roughly means that $\g(x|X)\,dx$ is proportional 
to the conditional probability for the existence of a particle 
in an infinitesimal volume $dx$ when the 
remaining configuration is $X$. Formally, it is not difficult to see 
that \rf{Papangelou} is equivalent to the statement that $P$
is absolutely continuous with respect to the intensity-$1$ Poisson 
point process $\pi_\L=\pi_\L^1$ in $\L$ with Radon--Nikodym density $g$ 
satisfying $g(X\cup\{x\})=\g(x|X)\,g(X)$, see e.g. Georgii and 
K\"uneth (1997). In particular, 
the Poisson process $\pi_\L^z$ of intensity $z>0$ on $\L$ has 
Papangelou intensity $\g(x|X)=z$.
\index{stochastic order!of point processes} 
\begin{prop}\label{Holl-Pres}
{\bf(Holley-Preston inequality) }Let $\L\subset\Rd$ be a bounded 
Borel 
set and $P,P'$ two probability measures on $\X_\L$ with Papangelou 
intensities $\g$ resp.\ $\g'$. 
Suppose $\g(x|X)\le \g'(x|X')$ whenever $X\subset X'$ and $x\notin 
X'\setminus X$. Then $P\preceq P'$. 
If this condition holds with $P'=P$ then $P$ has positive 
correlations.
\end{prop} 
Of course, the stochastic partial order $P\preceq P'$ is 
defined by means of the inclusion relation on $\X_\L$. 
Under additional technical assumptions the preceding proposition was 
first derived by Preston in 1975;
in the present form it is due to Georgii and 
K\"uneth (1997).
In the next section we will see how this result can be used to 
establish percolation in certain continuum random-cluster models, and 
thereby the existence of phase transitions in certain continuum 
particle systems.

\section{The continuum Ising model}
\index{Ising model!continuum|(}

The continuum Ising model is a model of point particles in $\Rd$ 
of two different 
types, plus and minus. Rather than of particles of different types, 
one may also think of particles with a ferromagnetic spin with two 
possible orientations. The latter would be suitable for modelling 
ferrofluids such as the Au-Co alloy, which have recently found some 
physical attention. Much of what follows can also be extended to 
systems with more than two types, but we stick here to the simplest
case. A configuration of particles is then described by a pair
$\xi=(X^+,X^-)$, where $X^+$ and $X^-$ are the configurations of 
plus- 
resp.\ minus-particles. The configuration space is thus $\O=\X^2$.

We assume that the particles interact via a 
repulsive interspecies pair potential of finite range, which is given 
by an even measurable function $\ph:\Rd\to[0,\infty]$ of bounded 
support. The Hamiltonian in a bounded Borel set $\L\subset\Rd$ 
of a configuration $\xi=(X^+,X^-)$ is thus given by
\be{Cont-Ising}
H_\L(\xi)=\sum_{x\in X^+, y\in X^-:\, \{x,y\}\cap\L\ne\emptyset}
\ph(x-y)\;.
\ee
In view of its analogy to \rf{Ising} this model is called the 
{\em continuum Ising model}. Setting $\ph(x-y)=\infty$ when 
$|x-y|\le 2r$ and $\ph(x-y)=0$ otherwise, we obtain the
classical {\em Widom--Rowlinson model} (1970) with 
a hard-core interspecies repulsion
\index{Widom--Rowlinson!continuum model}  
(which in spatial statistics is occasionally referred to as
the {\em penetrable spheres mixture model\/}). 
Of course, this case corresponds to the Widom--Rowlinson 
lattice gas considered above. 
Here we make only the much weaker assumption that $\ph$ is 
bounded away from zero on a neighborhood of the 
origin. That is, there exist constants $\d,r>0$ such that 
\be{pos}
\ph(x-y)\ge\d\mbox{ when }|x-y|\le 2r.
\ee

An interesting generalization of the Hamiltonian \rf{Cont-Ising} can 
be obtained by adding an interaction term which is independent of the 
types of the particles. In a ferrofluid model this would mean that in 
addition to the ferromagnetic interaction of particle spins there is 
also a molecular interaction which is spin-independent. Such 
an extension is considered in Georgii and H\"aggstr\"om (1996).

Given the Hamiltonian \rf{Cont-Ising}, the associated 
{\em Gibbs distribution in $\L$ 
with activity $z>0$ and boundary condition 
$\xi_\Lc=(X^+_\Lc,X^-_\Lc)\in\X_\Lc^2$\/} is defined by the formula
\be{cont-Gibbs}
 G_\L(d\xi_\L|\xi_\Lc)= Z_{\L|\xi_{\L^c}}^{-1} 
 \exp[- H_\L(\xi)] \; \pi^z_\L(dX^+_\L)\,\pi^z_\L(dX^-_\L)
\ee
which is completely analogous to \rf{Gibbs}. The corresponding set 
$\G=\G(z)$ 
of Gibbs measures is then defined as in Definition \ref{DLR}
(with $\L$ running through the bounded Borel sets in $\Rd$ instead of 
the finite subsets of $\Zd$).
\index{Gibbs measure!continuum}

In general, the existence of Gibbs measures in continuum models is 
not easy to establish. In the present case, however, it is simple:
Thinking of $\X_\L^2$ as the space of configurations on two disjoint 
copies $\L^+$ and $\L^-$ of $\L$, we see that
$G_\L(\,\cdot\,|\xi_\Lc)$ has the Papangelou intensity 
\be{Poisson-dom}
\g(x|X^+,X^-)= \left\{\ba{cl}
z\,\exp[-\sum_{y\in X^-}\ph(x-y)]&\mbox{if $x\in\L^+$}\\
z\,\exp[-\sum_{y\in X^+}\ph(x-y)]&\mbox{if $x\in\L^-$}\ea\right\} \le 
z\;.
\ee
Proposition \ref {Holl-Pres} therefore implies that
 $G_\L(\,\cdot\,|\xi_\Lc)\preceq \pi^z_\L\times \pi^z_\L$. Standard 
compactness theorems for point processes now show that for 
each $\xi\in\X^2$ the sequence $G_\L(\,\cdot\,|\xi_\Lc)$ has an 
accumulation 
point $P$ as $\L\uparrow\Rd$, and it is easy to see that $P\in\G$.

\subsection{Uniqueness and phase transition}

We will now show that the Gibbs measure is unique 
when $z$ is small, whereas a phase transition occurs when $z$ 
is large enough. Both results rely on percolation techniques.
As in the Widom--Rowlinson lattice gas, it remains open
whether there is a sharp activity threshold 
separating intervals of uniqueness and non-uniqueness.
\begin{prop}\label{unique}
For the continuum Ising model we have $\#\G(z)=1$ when $z$ is 
sufficiently small.   
\end{prop}
{\em Sketch proof: }Let $P,P'\in \G(z)$. We show that $P=P'$ when $z$ 
is small enough. Let $R$ be the range of $\ph$, i.e., 
$\ph(x)=0$ when $|x|\le R$, and divide $\Rd$ into cubic cells 
$\D(i)$, 
$i\in\Zd$, of linear size $R$. Let $p_c^*$ be the Bernoulli site 
percolation threshold of the graph with vertex set $\Zd$ and edges 
between all points having distance 1 {\em in the max-norm}. 
Consider the Poisson measure $Q^z=\pi^z\times\pi^z$ on the 
configuration space $\O=\X^2$. 

Let $\xi,\xi'$ be two independent realizations of $Q^z$, and
suppose $z$ is so small that $Q^z\times 
Q^z(N_i+N_i'\ge1)<p_c^*$, where $N_i$ and $N_i'$ are the numbers of 
particles (plus or minus) in  $\xi$ resp.\ $\xi'$.
Then for any finite union $\L$ of cells we have
$Q^z\times Q^z(\L\llr{\ge1}\infty)=0$, where 
$\{\L\llr{\ge1}\infty\}$ denotes the event that a cell in $\L$ 
belongs to an 
infinite connected set of cells $\D(i)$ containing at least one 
particle in either $\xi$ or $\xi'$.
Proposition \ref{Holl-Pres} together with \rf{Poisson-dom} imply that 
$P\times P'\preceq Q^z\times Q^z$. 
Hence $P\times P'(\L\llr{\ge1}\infty)=0$.
In other words, given two 
independent realizations $\xi$ and $\xi'$ of $P$ and $P'$ 
there exists a random corridor 
of width $R$ around $\L$ which is completely free of particles. In 
particular, this means that $\xi$ and $\xi'$ coincide on this random 
corridor.
By a spatial strong Markov property of Gibbs measures, it follows 
that $P$ and $P'$ coincide on the $\s$-algebra of events in $\L$.
As $\L$ can be chosen arbitrarily large, this proves the proposition.
$\Box$

\medskip\noindent
After this result on the absence of phase transition (following from 
the absence of some kind of percolation) we turn to the 
existence of phase transition. This will follow from the existence of 
percolation in a suitable random-cluster model. In analogy to 
Propositions \ref{prop:RC} and \ref{prop:WR-RC}, we will derive a
\index{random-cluster representation!continuum Ising model|(}
random-cluster representation of the Gibbs distribution
\be{cont-plus}
G_\L^+=\int \pi^z_\Lc(dY^+_\Lc)\;G_\L(\,\cdot\,| Y_\Lc^+,\emptyset)
\ee
with a Poisson boundary condition of plus-particles and no 
minus-particle off $\L$. 
Its random-cluster counterpart is the following probability measure
$\chi_\L$ on $\X\times\E$ describing random graphs $(Y,E)$ in $\Rd$:
\be{cont-RC}
\chi_{\L}(dY,dE) = Z_{\L|Y_\Lc}^{-1}\; 2^{k(Y,E)}\, 
\pi^z(dY)\,\mu^{p,\L}_{Y}(dE)\;.
\ee
In the above, $k(Y,E)$ is the number of clusters of the graph 
$(Y,E)$, 
$Z_{\L|Y_\Lc}=\int 2^{k(Y,E)}\,\pi^z_\L(dY_\L)$ normalizes the 
conditional probability of $\chi_\L$ given $Y_\Lc=Y\cap\Lc$ (so that
$Y_\Lc$ still has the Poisson distribution $\pi^z_\Lc$), and
$\mu^{p,\L}_{Y}$ is the probability measure on $\E$ for which the 
edges $e=\{x,y\}\subset Y$ are drawn independently with probability 
$p(x-y)=1-e^{-\ph(x-y)}$ if $e\not\subset Y_\Lc$, and probability $1$ 
otherwise. The probability measure $\chi_{\L}$ in \rf{cont-RC}
is called the {\em continuum random-cluster distribution in $\L$ with 
connection probability function $p$ and wired boundary condition}.
Note that (in contrast to \rf{RC} and \rf{WR-RC}) this distribution 
describes random configurations of both points and edges. 
In the Widom--Rowlinson case of a hard-core interspecies repulsion 
\index{Widom--Rowlinson!continuum model} the randomness of the 
edges disappears, and $\chi_{\L}$ describes a dependent Boolean 
percolation model which is the direct continuum analog of 
\rf{WR-RC}. The
{\em random-cluster representation of the continuum Ising model} now 
reads as follows.
\begin{prop}\label{prop:cont-RC} 
For any bounded box $\L$ in $\Rd$ there is the following 
correspondence between the Gibbs distribution $G_\L^+$ in 
\rf{cont-plus} for the continuum Ising model 
and the random-cluster distribution $\chi_\L$ in \rf{cont-RC}.

$(\,G_\L^+\leadsto\chi_\L\,) $ Take a particle configuration 
$\xi=(X^+,X^-)\in\O$ with distribution $G_\L^+$ and define a random 
graph $(Y,E)\in\X\times\E$ as follows: Let $Y=X^+\cup X^-$, and
independently for all $e=\{x,y\}\in E(Y)$ let $e\in E$ with 
probability 
\[
p_\L(e)=\left\{\ba{cl}
1-e^{-\ph(x-y)}&\mbox{\rm if $e\subset X^+$ or $e\subset X^-$, 
and $e\cap\L\ne\emptyset$,}\\ 
1&\mbox{\rm if $e\subset\Lc$,}\\
0&\mbox{\rm otherwise.} \ea\right.
\]
Then $(Y,E)$ has distribution $\chi_\L$.

$(\,\chi_\L\leadsto G_\L^+\,) $  Pick a random graph 
$(Y,E)\in\X\times\E$ according to $\chi_\L$. Define a particle 
configuration $\xi=(X^+,X^-)\in\O$ with $X^+\cup X^-=Y$ as follows:
For each finite cluster $C$ of $(Y,E)$ let $C\subset X^+$
or $C\subset X^-$ according to independent flips of a fair coin; 
the unique infinite cluster of $(Y,E)$ containing $Y_\Lc$ is included 
into $X^+$. Then $\xi$ has distribution $G_\L^+$.
\end{prop}
Just as in the lattice case, the random-cluster representation above 
gives the following key identity: for any finite box $\D\subset\L$,
\beal{cont-WR-percol}
\lefteqn{\int [\# X^+_\D -\# 
X^-_\D]\;G_\L^+(dX^+,dX^-)}\hspace*{8mm}\nonumber\\
&&=\int 
\#\{x\in Y_\D:x\lr Y_\Lc\}\;\chi_{\L}(dY,dE)\;;
\eeal
in the above, the notation $x\lr Y_\Lc$ means that $x$ is 
connected to a point of $Y_\Lc$ in the graph $(Y,E)$. 
\index{random-cluster representation!continuum Ising model|)}
In other words, the difference between the mean number of plus- and 
minus-particles in $\D$ corresponds to the percolation probability in 
$\chi_\L$. How can one check that the latter is positive for 
large $z$? The idea is again a stochastic comparison.

Let $\nu_\L=\chi_\L(\cdot \times\E)$ the point 
marginal of $\chi_\L$. Then
$\chi_\L(dY,dE)=\nu_\L(dY)\,\phi_{\L,Y}(dE)$ with an 
obvious analog $\phi_{\L,Y}$ of \rf{RC}. An application of 
Proposition \ref{Holl} shows that $\phi_{\L,Y}$ is stochastically 
larger than the Bernoulli edge measure $\mu_Y^{\tilde p}$ for which 
edges 
are drawn independently between points $x,y\in Y$ with probability 
$\tilde p(x-y)=(1-e^{-\d})/[(1-e^{-\d})+2\,e^{-\d}]$ when $|x-y|\le 
2r$, 
and with probability $0$ otherwise. Here, $\d$ and $r$ are as in 
assumption \rf{pos}. Moreover, $\nu_\L$ has the Papangelou intensity
\[
\g(x|Y)= z\int 2^{k(Y\cup\{x\},\,\cdot\,)}\, d\phi_{\L,Y\cup\{x\}}
\Big/\int 2^{k(Y,\,\cdot\,)}\, d\phi_{\L,Y}\;.
\]
To get a lower estimate for $\g(x|Y)$ one has to compare the effect 
on the number of 
clusters in $(Y,E)$ when a particle at $x$ and corresponding edges 
are added. In principle, this procedure could connect a large number
of distinct clusters lying close to $x$, so that 
$k(Y\cup\{x\},\,\cdot\,)$ was much smaller than $k(Y,\,\cdot\,)$.
However, one can show that this occurs 
only with small probability, so that $\g(x|Y)\ge \a z$ for some $\a>0$.
By Proposition \ref{Holl-Pres}, we can conclude that 
$\chi_\L$ is stochastically larger than the Poisson 
random-edge measure $P^{\a z,\tilde p}$ defined in \rf{rand-edge}. 
The right-hand side of \rf{cont-WR-percol} is therefore not smaller 
than $\theta(\a z,\tilde p;\Rd)$. Finally, since 
$G^+_\L\preceq \pi^z\times\pi^z$ by \rf{Poisson-dom}, the
Gibbs distributions $G^+_\L$ have 
a cluster point $P^+\in\G(z)$ satisfying
\[
\int [\# X^+_\D-\# X^-_\D]\;P^+(dX^+,dX^-)\geq 
\theta(\a z,\tilde p;\Rd)\;.
\]
By spatial averaging one can achieve that $P^+$ 
is in addition translation invariant. Together with 
Theorem \ref{cont-perc} this leads to the following theorem.
\begin{thm}\label{cont-phtr}
For the continuum Ising model on $\Rd$, $d\ge2$, 
with Hamiltonian \rf{Cont-Ising} and sufficiently large activity
$z$ there exist two translation invariant Gibbs measures $P^+$ and 
$P^-$ having a majority of plus- resp.\ minus-particles and related 
to each other by the plus-minus interchange.
\end{thm}
\index{phase transition!continuum Ising model} 
This result is due to Georgii and H\"aggstr\"om (1996). In the 
special case of the Widom-Rowlinson model 
\index{Widom--Rowlinson!continuum model} it has been derived 
independently in the same way by Chayes, Chayes, and 
Koteck\'y (1995). The first proof of phase transition in 
the Widom-Rowlinson model was found by Ruelle 
in 1971, and for a soft but strong repulsion by Lebowitz and Lieb in 
1972. Gruber and Griffiths (1986) used a
 direct comparison with the lattice Ising model in the case of 
a species-independent background hard core. 

As a matter of fact, one can make further use of stochastic 
monotonicity. (In contrast to the preceding theorem, this only works 
in the present case of two particle types.) 
Introduce a partial order `$\le$' on $\O=\X^2$ by writing
\be{anti-order}
(X^+,X^-)\le (Y^+,Y^-) \mbox{ when $X^+\subset Y^+$ and $X^-\supset Y^-$.}
\ee
A straightforward extension of Proposition \ref{Holl-Pres} then shows 
that the measures $G^+_\L$ in \rf {cont-plus} decrease stochastically 
relative to this 
order when $\L$ increases. (This can be also deduced from the 
couplings obtained by perfect simulation, see Section \ref{simul} 
below.) It follows that $P^+$ is in fact the limit 
of these measures, and is in particular translation invariant. 
Moreover, one can see that $P^+$ is stochastically maximal in 
$\G$ in this order.
This gives us the following counterpart to Corollary \ref{WR-plus-perc}.
\begin{cor}\label{cont-phtr-perc} 
For the continuum Ising model with any activity  $z>0$, a phase 
transition occurs  if and only if 
$$\int \hat P^+(dX^+,dX^-)\, \mu^p_{X^+}
(0\llr{+}\infty)>0\;;
$$
here $\hat P^+$ is the Palm measure of $P^+$, and
the relation $0\llr{+}\infty$ means that the origin 
belongs to an infinite cluster in the graph with vertex set $X^+$ and 
random edges drawn according to the probability function 
$p=1-e^{-\ph}$.
\end{cor} \index{percolation!continuum Ising model}

It is not known whether $P^+$ and $P^-$ are the only extremal 
elements of $\G(z)$ when $d=2$, as it is the case in the lattice 
Ising 
model. However, using a technique known in physics as the 
Mermin--Wagner theorem one can show the following.
\begin{thm} 
If $\ph$ is twice continuously differentiable then each $P\in\G(z)$ 
is translation invariant.
\end{thm}
A proof can be found in Georgii (1999). The existence of 
non-translation invariant Gibbs measures in dimensions $d\ge 3$ is an 
open problem.

\subsection{Thermodynamic aspects}
\label{thdyn}

Although we were able to take some advantage of stochastic comparison 
techniques in the continuum Ising model, the use of Proposition 
\ref{Holl-Pres} is much more limited than that of
its lattice analog. The reason is that its condition requires some 
kind of attractivity, which is in conflict with stability (preventing the 
existence of infinitely many particles in a bounded region).
This implies that a continuum Gibbs 
distribution $G_\L$ with a pair interaction cannot satisfy the 
conditions of
Theorem \ref{Holl-Pres} with $P=P'=G_\L$, which would imply that
$G_\L$ has positive correlations and is stochastically increasing 
with the activity $z$. Fortunately, this gap can be closed to 
some extent by the use of classical convexity techniques of 
Statistical Mechanics. These will allow us to conclude that 
at least the particle density of Gibbs 
measures is an increasing function of the activity $z$.

It should be noted that these ideas are standard in Statistical 
Physics; they are
included here because they might be less known among spatial 
statisticians, and because we need to check that the general 
principles really work in the model at hand. One should also note that 
this technique does not depend on the specific features of the 
model; in particular, 
it applies also to the gas of hard balls discussed in 
H.\ L\"owen's contribution to this volume.
\index{hard-core process}

Let us begin recalling the thermodynamic justification of Gibbs 
measures.
Let $\L\subset\Rd$ be a finite box. For any  translation invariant 
probability measure $P$ on $\O=\X^2$ consider the {\em entropy 
per volume} 
\[
\ss(P)=\lim_{|\L|\to\infty} |\L|^{-1}\, S(P_\L)\;.
\]
Here we write $P_\L$ for the restriction of $P$ to 
$\X_\L^2$, the set of particle configurations in $\L$, and
\[
S(P_\L)=\left\{\ba{cl}
-\int \log f \; dP_\L&\mbox{if $P_\L\ll \pi^1_\L\times\pi^1_\L$ 
with density $f$,}\\
-\infty&\mbox{otherwise}
\ea\right.
\]
for the entropy of $P_\L$ relative two the 
two-species Poisson process $\pi^1_\L\times\pi^1_\L$ on $\X_\L^2$ 
with intensity $1$. The notation $|\L|\to\infty$ means that $\L$ runs 
through a 
specified increasing sequence of cubic boxes with integer sidelength.
The existence of $\ss(P)$ is a multidimensional version of Shannon's 
theorem; see Georgii (1988) for the lattice case to which the present 
case can be reduced by identifying $\O$ with $(\X_C^2)^\Zd$ 
for a unit cube $C$.

Next consider the {\em interaction energy per volume}
\be{u}
\uu(P)=\int \hat P(dX^+,dX^-) 
\Big[1_{\{0\in X^+\}}\sum_{x\in X^-}\ph(x)
+ 1_{\{0\in X^-\}}\sum_{x\in X^+}\ph(x)\Big]\;
\ee
defined in terms of the Palm measure $\hat P$ of $P$. $\uu(P)$ 
can also be defined as a per-volume limit, cf.\ Georgii (1994), 
Section 3.
Also, consider the {\em particle density}
\[ \index{particle density}
\r(P)= \hat P(\O)=|\L|^{-1}\,\int[\# X^+_\L+\# X^-_\L]\, dP_\L
\]
of $P$; by translation invariance the last term does not depend on 
$\L$. The term $-\r(P)\,\log z$ is then equal to the chemical 
energy per volume.

Finally, consider the {\em  pressure} 
\be{p}
\pp(z)=-\min_P \Big[\uu(P)-\r(P)\,\log z -\ss(P)\Big]\;;
\ee
the minimum extends over all 
translation invariant probability measures $P$ on $\O$.
The large deviation techniques of Georgii (1994) show that 
$\pp(z)=\lim_{|\L|\to\infty} |\L|^{-1}\, 
\log Z_{\L|\xi_\Lc}$ for each $\xi\in\O$.
(The paper Georgii (1994) deals only with particles of a 
single type and superstable interaction, but the extension to the 
present case is straightforward because $\ph$ is nonnegative and has 
finite range.) 
The variational principle for Gibbs measures then reads as 
follows. \index{variational principle}
\begin{thm}\label{varprinc} 
Let $P$ be a translation invariant probability measure on $\O=\X^2$. 
Then $P\in\G(z)$ if and only if $\uu(P)-\r(P)\,\log z -\ss(P)$,
the free energy per volume, is equal to its minimum $-\pp(z)$.
\end{thm}
The ``only if'' part can be derived along the 
lines of Georgii (1994) and Proposition 7.7 of Georgii (1995). 
The ``if'' part 
follows from the analogous lattice result (see Section 15.4 of 
Georgii (1988)) by the identification of $\O$ and 
$(\X_C^2)^\Zd$ mentioned above.

What does the theorem tell us about the particle densities of Gibbs 
measures? Let us look at the pressure $\pp(z)$.
First, it follows straight from the definition \rf{p} 
that  $\pp(z)$ is a convex function of $\log z$. 
In other words, the function $\tilde{\pp}(t)=\pp(e^t)$ is convex.
Next, inserting $P=\pi^z\times\d_\emptyset$ into the right-hand 
side of \rf{p} we see that $\tilde{\pp}>-\infty$, and that the slope 
of $\tilde{\pp}$ at $t$ tends to infinity as $t\to\infty$.

Now, suppose $P\in\G(z)$ is translation invariant. The 
variational principle above then implies that the function 
$t\to (t-\log z)\rho(P) +\pp(z)$ is a tangent to $\tilde{\pp}$ at 
$\log z$. For, on the one hand we have
$$\uu(P)-\r(P)\,\log z -\ss(P)=-\pp(z)<\infty$$ 
and thus $\uu(P)-\ss(P)<\infty$, and on the other hand
\[
\uu(P)-\r(P)\,t -\ss(P)\geq-\tilde{\pp}(t) \quad\mbox{ for all $t$.}
\]
Inserting the former identity into the last inequality we get
the result. As a consequence, the particle density $\r(P)$ 
lies in the interval between the left and right derivative of 
$\tilde{\pp}$ at $\log z$. By convexity, these derivatives are 
increasing and almost everywhere identical. In fact, they are 
strictly increasing. For, if they were constant on some non-empty 
open interval $I$ then for each $t_0\in I$ and $P\in\G(e^{t_0})$ 
the function $t\to \r(P)\,t-\tilde{\pp}(t)$ would be constant on $I$, 
and thus by the variational principle $P\in\G(e^{t})$ for all 
$t\in I$. This is impossible because the conditional Gibbs 
distributions depend non-trivially on the activity.
We thus arrive at the following conclusion.
\begin{cor}\label{monotone} 
Let $0<z<z'$ and $P\in\G(z)$, $P'\in\G(z')$ be translation invariant.
Then $\r(P)< \r(P')$, and $\r(P)\to\infty$ as $z\to\infty$.
\end{cor} \index{particle density!monotonicity}
In the present two-species model it is natural to consider also 
the case when each particle species has its own activity, i.e., 
the plus-particles have activity $z^+$ and the minus-particles 
have activity $z^-$. It then follows in the same way that the pressure
$\pp(z^+,z^-)$ is a strictly convex function of the pair 
$(\log z^+,\log z^-)$, and therefore that
the density of plus-particles is a strictly increasing function 
of $z^+$, and the density of minus-particles is a strictly 
increasing function of $z^-$; these densities tend to infinity 
as $z^+$ resp.\ $z^-$ tends to infinity.
Moreover, Theorem \ref{cont-phtr} implies that $\pp(z^+,z^-)$ has 
a kink at
$(z,z)$ when $z$ is large enough; this  means that the convex 
function $t\to \pp(ze^t,ze^{-t})$ is not differentiable at $t=0$.

\subsection{Projection on plus-particles}

As the continuum Ising model is a two-species model, it is natural 
to ask what kind of system appears if we forget all minus-particles 
and only retain the plus-particles. The answer is that their 
distribution is again 
Gibbsian for a suitable Hamiltonian. This holds also in the case of 
different activities $z^+$ and $z^-$ of plus- and minus-particles, 
which is the natural context here.
To check this, take any box 
$\D\subset\Rd$ and let $\L\supset\D$ be so large that the distance of 
$\D$ from $\Lc$ exceeds the range $R$ of $\ph$.
Integrating over $X_\L^-$ in \rf{cont-Gibbs} and conditioning on 
$X^+_{\L\setminus\D}$ one finds
that the conditional distribution of $X^+_\D$ for given 
$X^+_{\L\setminus\D}$ under $G_\L(\,\cdot\,|\xi_\Lc)$ does not depend 
on $\L$ and $\xi_\Lc$ and has the Gibbsian form
\be{plus-Gibbs}
\p{G}_\D(dX^+_\D|X^+_{\D^c})= {\p{Z}_{\D|X^+_{\D^c}}}^{-1} 
\exp[-z^-\;\p{H}_\D(X^+)]  \; \pi^{z^+}_\D(dX^+_\D)
\ee
with the Hamiltonian
\be{plus-Ham}
\p{H}_\D(X^+)= \int_\D \Big( 1-\exp\Big[-\sum_{x\in 
X^+}\ph(x-y)\Big]\Big)\; dy \;.
\ee
Thus, writing $\G(z^+,z^-)$ for the set of all continuum-Ising Gibbs 
measures on $\O=\X^2$ with Hamiltonian \rf{Cont-Ising} and activities 
$z^+$ and $z^-$, 
and $\p{\G}(z^+,z^-)$ for the set of all Gibbs measures on $\X$ with 
conditional distributions \rf{plus-Gibbs}, one obtains the following 
corollary.
\begin{cor}\label{plus-part} 
Let $P\in\G(z^+,z^-)$ and $\p{P}$ be the 
distribution of the configuration of plus-particles. Then
$\p{P}\in\p{\G}(z^+,z^-)$. In particular, 
$\#\p{\G}(z,z)>1$ when $z$ is large enough.
\end{cor}
The last statement follows from Theorem \ref{cont-phtr}.

In the Widom--Rowlinson case when $\ph=\infty\,1_{\{|\,\cdot\,|\le 
2r\}}$, the relationship between $\G(z^+,z^-)$ and $\p{\G}(z^+,z^-)$ 
has already \index{Widom--Rowlinson!continuum model}
been observed in the original paper by Widom and Rowlinson (1970). 
The plus-Hamiltonian \rf{plus-Ham} then takes the simple form
\[
\p{H}_\D(X^+)= |\D\cap \bigcup_{x\in X^+} B_{2r}(x)|\;
\]
that is, $\p{H}_\D(X^+)$ is the volume in $\D$ of the Boolean model 
with 
radius $2r$ induced by $X^+$. In this case, the model with 
distribution
$\p{G}_\D(dX^+_\D|X^+_{\D^c})$ was reinvented by Baddeley and van
Lieshout (1995). Having the two-dimensional case in mind, they 
coined the suggestive term {\em area-interaction process}.
\index{area-interaction process} From here 
one can go one step further to Hamiltonians which use not only the 
volume but also the other Minkowski functionals. This has been 
initiated by Likos et al. (1995) and Mecke (1996); see also Mecke's 
contribution to this volume. 

One particularly nice feature of the area-interaction process and its 
generalization \rf{plus-Ham} is that it seems to be the only known 
(non-Poisson) model to which Proposition \ref{Holl-Pres} can be 
applied for 
establishing positive correlations of increasing functions. This 
attractiveness property makes the model quite attractive for 
statistical modelling. (By way of contrast, repulsive point systems 
can be modelled quite easily, for example by a nonnegative pair 
interaction.)

However, some caution is necessary 
due to the phase transition when $z^+=z^-=z$ is large: 
The typical configurations of $\p{G}_\D(\,\cdot\,|\,\emptyset\,)$ for a 
large finite window $\D$ then can be typical for either phase, $\p{P}^+$ 
or $\p{P}^-$, and thus can have different particle densities. Due to 
finite size effects, this phenomenon already appears when $z^+$ and 
$z^-$ are sufficiently close to each other. So, the spatial 
statistician should be aware of 
such an instability of observations and should examine whether this 
is realistic or not in the situation to be modelled.

Finally, one can use Proposition \ref{Holl-Pres} to show that the
Gibbs measures $\p{P}\in\p{G}(z^+,z^-)$ are stochastically increasing
in $z^+$ and decreasing in $z^-$. In particular, the density 
of plus-particles for any $P\in\G(z^+,z^-)$ increases when 
$z^+$ increases or $z^-$ decreases, as can also be seen using
the partial order \rf{anti-order}. 
\index{particle density!monotonicity} The monotonicity
results in the last paragraph of Section \ref{thdyn} thus follow
also from stochastic comparison techniques, but
Corollary \ref{monotone} cannot be derived in this way.

\subsection{Simulation}
\label{simul}

There are various reasons for performing Monte--Carlo simulations of 
physical or statistical systems, as discussed in a number of 
other contributions to this volume. In the present context, the 
primary reason is to sharpen the intuition on the system's behavior, 
so that one can see which properties can be expected to hold. This 
can lead to conjectures which then hopefully can be 
checked rigorously.

Here we will show briefly how one can obtain simulation pictures of 
the continuum Ising model. We start with a continuum Gibbs sampler 
which is suggested by Proposition \ref{prop:cont-RC}; in the 
Widom--Rowlinson case it has been proposed by H\"aggstr\"om, van Lieshout 
and M\o ller (1997). \index{Widom--Rowlinson!continuum model}

Consider a fixed window $\L\subset\Rd$ and the Gibbs distribution 
$G_\L^{\mbox{\scriptsize free}}=G_\L(\,\cdot\,|\emptyset,\emptyset)$ 
for the continuum 
Ising model in $\L$ with activity $z$ and free (i.e., empty) 
boundary condition off $\L$. 
We define a random map $F:\X_\L\to\X_\L$ by the following algorithm:
\bit
\item take an input configuration $X\in\X_\L$,
\item select a Poisson configuration $Y\in\X_\L$ with distribution 
$\pi_\L^z$,
\item define a random edge configuration $E\subset \{\{x,y\}:x\in X, 
\,y\in Y\}$ by independently drawing an edge from $x\in X$ to $y\in Y$
with probability $p(x-y)=1-e^{-\ph(x-y)}$,
\item set $F(X)=\Big\{y\in Y:\{x,y\}\notin E \mbox{ for all }x\in X 
\Big\}$.
\eit
That is, $F(X)$ is a random thinning of $\pi^z_\L$ obtained by removing
all points which are connected to $X$ by a random edge. Its distribution
is nothing other than the
Poisson point process $\pi^{z,\ph}_{\L|X}$ on $\L$ 
with inhomogeneous 
intensity measure $\rho^\ph_X(dy)=z\,1_\L(y)\, \exp\left[-\sum_{x\in 
X}\ph(y-x)\right]\}\,dy$\,. (Of course, this could also be achieved by setting
$F(X)=\{y\in Y:U_y\leq \exp\left[-\sum_{x\in 
X}\ph(y-x)\right]\}$ for independent $U(0,1)$-variables $U_y$, $y\in Y$.
Although this was simpler in the case of the MCMC below, it would
considerably increase the running time of the perfect algorithm of
Theorem \ref{perfect}, as will be explained there.)
Now the point is that $\pi^{z,\ph}_{\L|X^+}$ is the conditional 
distribution of $X^-$ given $X^+$ relative to the Gibbs distribution
$G_\L^{\mbox{\scriptsize free}}$, and
similarly with $+$ and $-$ interchanged. So, if 
$\xi=(X^+,X^-)$ has distribution $G_\L^{\mbox{\scriptsize free}}$ and 
$F^+,F^-$ are 
independent realizations of $F$ then $(F^+(X^-), F^-\circ F^+(X^-))$ 
has again distribution $G_\L^{\mbox{\scriptsize free}}$. This 
observation gives rise to 
the following {\em Markov chain Monte Carlo algorithm} 
(MCMC). \index{Markov chain Monte Carlo!continuum Ising model}
\begin{prop}\label{MCMC} 
Let $F^+_n,F^-_n$, $n\ge 0$, be independent realizations of $F$, and 
$X^+_0\in\X_\L$ any initial configuration. Define recursively
\[
X^-_0=F^-_0(X^+_0), \quad X^+_n=F^+_n(X^-_{n-1}),\; 
X^-_n=F^-_n(X^+_n) \mbox{ for $n\ge1$.}
\]
Then the distribution of $(X^+_n,X^-_n)$ converges to 
$G_\L^{\mbox{\scriptsize free}}$ 
in total variation norm at a geometric rate.
\end{prop}
{\em Proof. }It suffices to observe that $F\equiv\emptyset$ with 
probability $\d=e^{-z\,|\L|}$. This shows that for any two 
configurations $X,X'\in\X_\L$ and any $A\subset\X_\L$
$$|\mbox{Prob}(F(X)\in A)-\mbox{Prob}(F(X')\in A)|\leq 1-\d\;.$$
So, if one looks at the process $(X^+_n,X^-_n)$ for two different 
starting configurations then each application of $F$ reduces the total 
variation distance by a factor of $1-\d$. $\Box$

\medskip\noindent 
A nice property of the random mapping $F$ is its 
monotonicity: if $X\subset X'$ then $F(X)\supset F(X')$ almost 
surely. This allows to modify the preceding algorithm to obtain {\em 
perfect simulation} in the spirit of Propp and Wilson, as 
described in the contribution of E.\ Th\"onnes to this 
volume. According to \rf{Poisson-dom}, 
$G_\L^{\mbox{\scriptsize free}}$ is stochastically dominated by 
independent Poisson processes of plus- 
and minus-particles. \index{perfect simulation!continuum Ising model|(}
So one can use the idea of {\em dominated 
perfect simulation} in her terminology. We describe here only the 
algorithm and refer to her contribution for more details.

Roughly speaking, the perfect algorithm consists of repeated simultaneous 
runs of the preceding MCMC, starting
from two particular initial conditions at some time
$N_k<0$ until time $0$. The two initial conditions are chosen extremal 
relative to the ordering \rf{anti-order}, namely with no initial 
plus-particle (the minimal case), and with a Poisson crowd of intensity $z$ 
of plus-particles (which is maximal by stochastic domination). Since 
the same realizations of $F$ are used in both cases, the two parallel 
MCMC's have a positive chance of coalescing during the time interval from
$N_k$ to $0$. If this occurs, one stops. Otherwise one performs 
a further run which starts at some time $N_{k+1}<N_k$.
\begin{thm}\label{perfect} 
Let $F^+_n,F^-_n$, $n\le 0$, be independent realizations of $F$, and
$(N_k)_{k\ge1}$ a strictly decreasing sequence of negative run
starting times. For each run indexed by $k\ge1$ let
\[
\Phi_k=F^+_0\circ F^-_{-1}\circ F^+_{-1}\circ\ldots\circ 
F^-_{N_k+1}\circ F^+_{N_k+1}\circ F^-_{N_k}
\]
be the random mapping corresponding to the MCMC of 
Proposition 
\ref{MCMC} for the time interval from $N_k$ to $0$, and
consider the processes 
$$X_{k,\min}^+=\Phi_k(\emptyset), \quad 
X_{k,\max}^+=\Phi_k\circ F^+_{N_k}(\emptyset)\;.$$
Then there exists a smallest (random) $K<\infty$ such that 
$X_{K,\min}^+=X_{K,\max}^+$, and the random particle configuration
$\xi_K=(X_{K,\min}^+, F^-_0(X_{K,\min}^+))$ has distribution 
$G_\L^{\mbox{\scriptsize free}}$.
\end{thm}
Since the random mapping $F$ can be simulated by simple standard 
procedures, the implementation of the preceding 
algorithm is quite easy; a Macintosh application can be found at
{\tt http://www.mathematik.uni-muenchen.de/\~{}georgii/CIsing.html }.
The 
main task is to store the random edge configuration $E$ in each 
application of $F$ during a time interval $\{N_k,\ldots,N_{k-1}-1\}$ 
for use in the later runs (which should be done in a file on the hard 
disk when $z|\L|$ is large). As a matter of fact, once the set $E$
is determined one can forget the positions of the particles of $Y$
and only keep their indices. In this sense, $E$ contains all 
essential information of the mapping $F$. As a consequence, 
knowing $E$ one needs almost no time to apply the same realization 
of $F$ in later runs. This is not the case for the alternative
definition of $F$ mentioned above.

However, there are some difficulties coming from the phase 
transition of the model. Running the 
perfect algorithm for small $z$ is fine and raises no problem. 
But if $z$ is large then the algorithm requires a considerably 
longer time to terminate. This is because for each 
run $k$ the distribution of 
$\xi_{k,\max}=(X_{k,\max}^+, F^-_0(X_{k,\max}^+))$ will be close to 
$P^+$ and thus 
show a large crowd of plus-particles giving the minus-particles 
only a minimal chance to spread out. Likewise, the distribution of 
$\xi_{k,\min}=(X_{k,\min}^+, F^-_0(X_{k,\min}^+))$ will be close to 
$P^-$, so that 
the minus-particles are in the great majority. The bottleneck between 
these two types of configurations is so small that $K$ will
typically be much too large for practical purposes, at least for 
windows $\L$ of satisfactory size. In order to reduce 
this difficulty, one should not simulate the 
Gibbs distribution $G_\L^{\mbox{\scriptsize free}}$ with free 
boundary condition (as we have done above for simplicity), as this
distributes most of its mass on two quite opposed events.
Rather one should simulate one of the phases, say $P^+$.
For a finite window,
this can be achieved by imposing a random boundary condition of 
Poisson plus-particles outside $\L$ as in 
\rf{cont-plus}. Such a boundary condition helps $X_{k,\max}^+$ and 
$X_{k,\min}^+$ quite a lot to coalesce within reasonable time. (If
one is willing to accept long running times, one should impose 
periodic boundary conditions to reduce the finite-size effects.)

\begin{center}
{\sl This electronic version does not contain the figures,\\ which can
be found on the website mentioned above.}
\end{center}

The pictures shown are obtained in this way. The underlying interaction 
potential is $\ph(x)=3\,(1-|x|)^2$ for $|x|\le 1$, $\ph(x)=0$ otherwise.
The size of the window is $20\times 20$. Outside of the window there 
is an invisible boundary condition of white Poisson particles.
The activities are $z=4.0$ (top) and $z=4.5$ (bottom), and the 
corresponding coalescence times are $-N_K=300$ resp.\ $400$. 

In the subcritical case $z=4$, the particles in the bulk do not feel the 
boundary condition: there is no 
dominance of white over black. In the supercritical case $z=4.5$
however, the influence of the white boundary condition is strong 
enough to dominate the whole window, and the phase transition becomes manifest.
This is nicely illustrated by
the random-cluster representation,
\index{random-cluster representation!continuum Ising model} which according to 
Proposition \ref{prop:cont-RC} is obtained from the point 
configuration $(X^+,X^-)$ by adding random edges within $X^+$ and 
$X^-$ separately. Here one sees that in the subcritical case the influence 
of the plus-boundary condition is only felt by a the particles near the 
boundary, while in the supercritical case the global behavior is 
dominated by a macroscopic cluster reaching 
from the boundary far into the interior of $\L$. This visualizes 
the equivalence of phase transition and percolation derived in 
Corollary \ref{cont-phtr-perc}. 

To conclude we mention two other algorithms. First, there is another 
perfect algorithm using a rejection scheme due to Fill, which has 
been studied in detail by Th\"onnes (1999) in the 
Widom-Rowlinson case; \index{Widom--Rowlinson!continuum model}
its extension to the present case is straightforward. 
A further possibility, which is particularly useful in the 
supercritical case, is to use a continuum analog of the 
Swendsen--Wang algorithm, Swendsen and Wang (1987). 
\index{Swendsen--Wang algorithm}
In its classical version for the lattice 
Ising model, this algorithm consists in alternating applications of 
the two procedures in Proposition \ref{prop:RC}. In its continuum version, 
one can again alternate between the two procedures of 
Proposition \ref{prop:cont-RC}, but one has to combine this with 
applications of the 
random mapping $F$ in order to obtain a resampling of particle 
positions. Unfortunately, this algorithm does not seem to admit a 
perfect version because of its lack of monotonicity, but it has the 
advantage of working also for the many-species extension of the model.
\index{Ising model!continuum|)}

\section*{References}
\small

\begin{list}{}{\setlength{\leftmargin}{0mm}}

\item Baddeley, A.J. and van Lieshout, M.N.M. (1995) 
Area-interaction point processes, {\sl Ann.\ Inst.\ Statist.\ 
Math.} {\bf 46}, 601--619.

\item Chayes, J.T., Chayes, L. and Koteck\'y, R. (1995)
The analysis of the Widom-Rowlinson model by stochastic geometric 
methods, {\sl Commun. Math. Phys.} {\bf 172}, 551--569.

\item Edwards, R.G. and Sokal, A.D. (1988) Generalization of the
Fortuin--\-Kasteleyn--\-Swendsen--\-Wang representation and Monte Carlo
algorithm, {\sl Phys. Rev.} D {\bf 38}, 2009--2012.

\item Fortuin, C.M. and Kasteleyn, P.W. (1972)
On the random-cluster model. I. Introduction and relation to other models,
{\sl Physica} {\bf 57}, 536--564.

\item Georgii, H.-O. (1988) {\sl Gibbs Measures and Phase
Transitions}, de Gruyter, Berlin New York.

\item  Georgii, H.-O. (1994) Large Deviations and the 
Equivalence of
Ensembles for Gibbsian Particle Systems with Superstable Interaction. 
{\sl Probab.\ Th.\ Rel.\ Fields} {\bf 99}, 171--195. 

\item Georgii, H.-O. (1995) The Equivalence of 
Ensembles for Classical Systems of Particles. {\sl J.\ Statist.\ 
Phys.} {\bf 80}, 1341--1378.  

\item Georgii, H.-O. (1999) 
Translation invariance and continuous symmetries in two-dimensional 
continuum systems, in: S. Miracle-Sole, J. Ruiz, V. Zagrebnov (eds.), 
{\sl Mathematical results in Statistical Mechanics}, Singapore etc., 
World Scientific, pp. 53--69. 

\item Georgii, H.-O. and H\"aggstr\"om, O. (1996)
Phase transition in continuum Potts models,
{\sl Commun. Math. Phys.} {\bf 181}, 507--528.

\item Georgii, H.-O., H\"aggstr\"om, O. and Maes, C. (1999)
The random geometry of equilibrium phases, in: Domb  and J.L. 
Lebowitz (eds.), {\sl Critical phenomena}, Academic Press, forthcoming.

\item Georgii, H.-O. and Higuchi, Y. (1999) Percolation 
and number of phases in the 2D Ising model, submitted to {\sl J.\ 
Math.\ Phys.}

\item Georgii, H.-O. and K\"uneth, T. (1997) Stochastic 
comparison of
point random fields, {\sl J. Appl. Probab.} {\bf 34}, 868--881.

\item Grimmett, G.R. (1999) {\sl Percolation}, 2.\ 
ed., Springer, New York.

\item Gruber, Ch. and Griffiths, R.B. (1986) Phase transition 
in a ferromagnetic fluid, {\sl Physica A} {\bf 138}, 220--230.

\item H\"aggstr\"om, O., van Lieshout, M.N.M. and M\o ller, 
J. (1997) Characterization results and Markov chain Monte Carlo 
algorithms including
exact simulation for some spatial point processes, {\sl Bernoulli}, to
appear.

\item Likos, C.N, Mecke, K.R. and Wagner, H. (1995) Statistical 
morphology of random interfaces in microemulsions, 
{\sl J. Chem. Phys.} {\bf 102}, 9350--9361.

\item Mecke, K.R. (1996) A morphological model for complex fluids, 
{\sl J. Phys.: Condens. Matter} {\bf 8}, 9663--9668.

\item Meester, R. and Roy, R. (1996) {\sl Continuum 
Percolation}, Cambridge University Press.

\item Penrose, M.D. (1991): On a continuum percolation model,
{\sl Adv. Appl. Probab.} {\bf 23}, 536--556.

\item Swendsen, R.H. and Wang, J.-S. (1987)
Nonuniversal critical dynamics in Monte Carlo simulations.
{\sl Phys. Rev. Lett.} {\bf 58}, 86--88.

\item Th\"onnes, E. (1999): Perfect simulation of some 
point processes for the impatient user,
{\sl Adv. Appl. Probab.} {\bf 31}, 69--87.

\item Widom, B. and Rowlinson, J.S. (1970) New model for the 
study of
liquid-vapor phase transition, {\sl J. Chem. Phys.} {\bf 52}, 
1670--1684.

\end{list}

%\printindex

\end{document}